\documentclass{cjour}
\usepackage{amssymb,psfig,url}

\def\Frac#1#2{\frac{\displaystyle{#1}}{\displaystyle{#2}}}
\begin{document}
 \authorrunninghead{A. Gil, W. Koepf $\&$ J. Segura}    %
 \titlerunninghead{Computation of zeros of hypergeometric functions}    %
\title{Numerical algorithms for the real zeros of hypergeometric functions}
\author{Amparo Gil$^{a}$, Wolfram Koepf$^{b}$ and Javier Segura$^{c}$}
\affil{$^a$Depto. de Matem\'aticas, Universidad Aut\'onoma de Madrid,
28049-Madrid, Spain \\
 E-mail: amparo.gil@uam.es\\
$^b$ Universit\"at Kassel, FB 17 Mathematik-Informatik, 34132-Kassel, Germany\\
  E-mail: koepf@mathematik.uni-kassel.de\\
$^c$ Depto. de Matem\'aticas, Estad\'{\i}stica y Computaci\'on.
Universidad de Cantabria. 39005-Santander,
Spain\\
  E-mail: javier.segura@unican.es   \\
}

\abstract{ Algorithms for the computation of the real zeros of
hypergeometric functions which are solutions of second order ODEs
are described. The algorithms are based on global fixed point
iterations which apply to families of functions satisfying first
order linear difference differential equations with continuous
coefficients. In order to compute the zeros of arbitrary solutions
of the hypergeometric equations, we have at our disposal several
different sets of difference differential equations (DDE). We
analyze the behavior of these different sets regarding the rate of
convergence of the associated fixed point iteration. It is shown
how combinations of different sets of DDEs, depending on the range
of parameters and the dependent variable, is able to produce
efficient methods for the computation of zeros with a fairly
uniform convergence rate for each zero. 

}



\begin{article}

\section{Introduction}

The zeros of hypergeometric functions are quantities which appear
in a vast number of physical and mathematical applications. For
example, the zeros of classical orthogonal polynomials (OP) are
the nodes of Gaussian quadrature; classical OP (Hermite, Laguerre
and Jacobi polynomials) are particular cases of hypergeometric
functions. Also, the zeros of Bessel functions and their
derivatives appear in many physical applications and there exists
a variety of methods of software for computing these zeros.

However, an efficient algorithm which can be applied to the computation
of all the zeros of any hypergeometric function in any real interval (not
containing a singular point of the defining ODE) is still missing.

In \cite{Seg2, Seg3} methods were introduced which are
capable of performing this task for hypergeometric functions which are solutions
of a second order ODE; an explicit Maple algorithm was presented in \cite{Gil}.
The starting point of the methods is the construction of a first order
system of differential equations
\begin{equation}
\begin{array}{ll}
y' (x)=\alpha (x) y(x) +\delta (x) w (x)\\
w' (x)=\beta (x) w(x) +\gamma (x) y (x),
\end{array}
\end{equation}
\noindent with continuous coefficients $\alpha (x)$, $\beta (x)$,
$\gamma (x)$ and $\delta (x)$ in the interval of interest,
relating our problem function $y(x)$ with a contrast function
$w(x)$, whose zeros are interlaced with those of $y(x)$.
Typically, the contrast function $w(x)$ satisfies a second order
ODE similar to the second order ODE satisfied by the problem
function.

Given a hypergeometric function $y(x)$ there are several known options to choose as contrast
functions $w(x)$. As an example, considering a Jacobi polynomial
\begin{equation}
y(x)=P_n^{(\alpha ,\beta)}(x)=\Frac{(\alpha+1)_{n}}{n!}\,_2\mbox{F}_{1} (-n,n+\alpha+\beta+1;
\alpha+1;(1-x)/2)
\end{equation}
\noindent we could take as contrast function $w_{\mbox{\tiny{OP}}} (x)=P_{n-1}^{(\alpha ,\beta)}(x)$
but also $w_{\mbox{\tiny{D}}}(x)=\Frac{d}{dx}P^{(\alpha ,\beta)}_n(x)$ is a possible choice. Both contrast
functions are $_2\mbox{F}_1 (a,b;c;x)$ hypergeometric functions with parameters $a,b,c$ differing
 by integer numbers from the parameters of the problem function $y(x)$.

When the contrast function in the previous example is $w_{\mbox{\tiny{OP}}}
(x)=P_{n-1}^{(\alpha ,\beta)}(x)$
the first order differential system is related to the three term recurrence relation for Jacobi
polynomials. It may seem that this is a natural differential system to consider. However, it was
numerically observed that the fixed point method which can be obtained from this differential
system becomes relatively slow for the zeros of $P_n^{(\alpha,\beta)}$ close to $\pm 1$ \cite{Gil}. Because
the extreme zeros approach to $\pm 1$ as $n\rightarrow +\infty$, the efficiency for the computation
of such zeros decreases as the order increases. Similar problems arise, for example, when
$\alpha \rightarrow -1^+$ or $\beta \rightarrow -1^+$. In fact, the number of iterations required
 to compute the extreme zeros tend to infinity in these limits. Similar problems take place for
 Laguerre polynomials $L_{n}^{\alpha }(x)$ for the smallest (positive) zero. Fortunately we will
 later show how the selection of $w_{\mbox{\tiny{D}}}(x)$ as contrast function
gives a much better asymptotic behavior for the resulting fixed
point iteration for the extreme Jacobi
 zeros. For the Laguerre case, a similar solution is possible.

These two examples illustrate the need to analyze the convergence
of the resulting fixed point iteration for the different available
contrast functions. Although for any adequate contrast function
(satisfying the necessary conditions \cite{Seg2,Seg3}) the
resulting fixed point method is quadratically convergent, the
non-local behavior of the method and the corresponding estimation
of first guess values for the zeros may result in disaster for
certain contrast
 functions in some limits. As a result of this study, we will
obtain explicit methods for the computation of the real zeros of
hypergeometric functions with a good asymptotic behavior and a
fairly uniform convergence rate in the whole range of parameters.

\section{Theoretical background}
\label{repaso}

Let us now briefly outline the main ingredients of the numerical method. For more details
we refer to \cite{Seg3,Seg2}.
It was shown in \cite{Seg3,Seg2} that, given a family of functions $\{y_k^{(1)}, y_k^{(2)}\}$,
depending on one parameter $k$, which are independent solutions of second order ODEs
\begin{equation}
y''_k+B_k(x)y'_k+A_k(x)y_k=0\,,k=n,n-1
\label{ODE}
\end{equation}
and satisfy relations of the type:
\begin{equation}
\begin{array}{rl}
y'_n=& a_n(x)y_n+d_n(x)y_{n-1}\\
y'_{n-1}= &b_n(x)y_{n-1}+e_n(x)y_n
\label{ddr}
\end{array}
\end{equation}
the coefficients $a_n(x)$, $b_n(x)$, $d_n(x)$, $e_n(x)$, $B_k (x)$
and $A_k(x)$ being continuous and $d_n e_n<0$ in a given interval
$[x_1, x_2]$, fixed point methods (Eq.\ (\ref{iter})) can be built
to compute all the zeros of the solutions of (\ref{ODE}) inside
this interval. These difference-differential equations (\ref{ddr})
are called general because they are satisfied by a basis of
solutions $\{y_k^{(1)}, y_k^{(2)}\}$. The fact that the DDEs are
general and with continuous coefficients in an interval $I$
implies \cite{Seg2} that $d_n e_n \neq 0$ in this interval.
Conversely, given $\{y_n^{(i)}, y_{n-1}^{(i)}\}$, $i=1,2$
independent solutions of the system (\ref{ddr}) and $d_n e_n\neq
0$, then $\{y_k^{(1)}, y_k^{(2)}\}$, $k=n, n-1$, are independent
solutions of the ODEs (\ref{ODE}). The method can then be applied
to compute the zeros of any solution of such ODEs.

It was shown that the ratios $H_i (z)$ ($i=\pm 1$):
\begin{equation}
\begin{array}{l}
H_i(z)=-i\,{\rm sign\:}(d_{n_i})K_{n_i}\Frac{y_n(x(z))}{y_{n_{i}}(x(z))}\\
K_{n_i}=\left(-\Frac{d_{n_i}}{e_{n_i}}\right)^{i/2}\,,\,
z(x)=\int \sqrt{-d_{n_i} e_{n_i}} dx
\end{array}
\label{FixP}
\end{equation}
where $n_{+1}=n+1\,,\,n_{-1}=n$, satisfy the first order equations
\begin{equation}
\dot{H}_i (z)=1+H_i (z)^2-2\eta_i (x(z)) H_i (z)
\label{ratH}
\end{equation}
\noindent where
$$
\eta_i (x)=i\Frac{1}{\sqrt{-d_{n_i}e_{n_i}}}\left(a_{n_i}-b_{n_i}+\frac12
\left(\Frac{e'_{n_i}}{e_{n_i}}-\Frac{d'_{n_i}}{d_{n_i}}\right)\right)
$$
and the dot means derivative with respect to $z$ while the prime
is the derivative with respect to $x$. \noindent Using Eq.\
(\ref{ratH}), one can show that
\begin{equation}
   T_i(z)=z - \arctan (H_i (z))
\label{iter}
\end{equation}
are globally convergent fixed point iterations (FPI): given a
value $z_0$ between two con\-se\-cu\-tive zeros
($z_{n_i}^{1},z_{n_i}^{2}$) of $y_{n_i}(x(z))$ (consecutive
singularities of $H_i (z)$),
 the iteration of (\ref{iter}) converges to $z_n\in (z_{n_i}^{1},z_{n_i}^{2})$, where $x_n=x(z_n)$
is a zero of $y_n (x)$.

Global bounds for the distance were provided which lead to iteration steps
that can be used to compute new starting values for obtaining all the zeros inside a given
interval.

It was shown that, in intervals where $\eta_i$ does not change sign,
 either
\begin{equation}
|z(x_{n_i}^{1})-z_n|<\pi /2 \mbox{ and  } |z(x_{n_i}^{2})-z_n|>\pi /2\,\,\,(\eta_i>0)
\label{etamas}
\end{equation}
\noindent
or
\begin{equation}
|z(x_{n_i}^{1})-z_n|>\pi /2 \mbox{ and  } |z(x_{n_i}^{2})-z_n|<\pi /2\,\,\, (\eta_i<0).
\label{etamenos}
\end{equation}
\noindent In this
way $\pi /2$ is the choice for the iteration step when the second situation takes place ($\eta_i<0$); this means
that if $y_{n_i}(x(z))$ has at least a zero larger than $z_n$, then
\begin{equation}
\lim_{j\rightarrow\infty}T^{(j)}(z_n+\pi/2 )
\label{nonimp}
\end{equation}
is the smallest zero larger than $z_n$. Similarly, when $\eta_i>0$, the iteration step will be $-\pi/2$ instead
of $\pi/2$ (backward sweep). When $\eta_i$ changes sign, forward and backward schemes can be combined \cite{Seg3}.

The functions $H_i (z)$ can be written as a ratio of functions $H_i (z)=\tilde{y}_{n}(z)/\tilde{y}_{n_i}(z)$, where
$\tilde{y}_k(z)=\lambda_k (z) y_k (x(z))$ and $\lambda_k$ have no zeros,
in such a way that the functions $\tilde{y}_n (z)$ and
$\tilde{y}_{n_i}(z)$ (with the same
zeros as $y_n (x(z))$ and $y_{n_i} (x(z))$ respectively) satisfy
 second order ODEs in normal form:
\begin{equation}
\begin{array}{ll}
\Frac{d^2 \tilde{y}_n}{dz^2}+\tilde{\mbox{A}}_n \tilde{y}_n=0\,, &
\Frac{d^2 \tilde{y}_{n_i}}{dz^2}+\tilde{\mbox{A}}_{n_i} \tilde{y}_{n_i}=0\,,  \\
& \\
\tilde{\mbox{A}}_n (z)=1+\dot{\eta}_i -\eta_i^2\,, & \tilde{\mbox{A}}_{n_i} (z)=1-\dot{\eta}_i -\eta_i^2\,.
\end{array}
\label{ODEsnor}
\end{equation}

Finally, we recall that using monotony conditions of $A_n (z)$ the
iteration steps $\pm \pi /2$ (Eq.\ (\ref{nonimp})) can be improved
according to Theorem 2.4 of \cite{Seg3}:

\begin{theorem}
If $z_{-1}<z_0<z_1$ are three consecutive zeros of $y_n (x(z))$ and $\eta_i (z)\dot{\tilde{\mbox{A}}}_n (z)>0$ in
$(z_{-1}, z_{1})$ then $z_{j}=\lim_{n\rightarrow\infty}T^{(n)}(z_0+\Delta z_0)$ where $\Delta z_0=z_0-z_{-j}$,
$j=\mbox{sign}(\eta)$. The convergence is monotonic.
\label{T2.4}
\end{theorem}

\subsection{Oscillatory conditions}

We are interested in computing zeros of oscillatory solutions of second order ODEs and, in particular, on building
algorithms for the computation of the zeros of the hypergeometric functions.
If a second order differential equation has a given number of singular regular points,
 we divide the real axis in subintervals determined by the singularities and search for
the zeros in each of these subintervals. We only apply the algorithms if it is not disregarded
 that the function can have two zeros at least
in the subinterval under consideration.

 We consider that an ODE has
oscillatory solutions in one of these subintervals if it has
solutions with at least two zeros in this subinterval; otherwise,
if all the solutions have one zero at most we will call these
zeros isolated zeros. The fixed point methods (FPMs) before
described deal with the zeros of any function satisfying a given
differential equation, no matter what the initial conditions are
on this function. Isolated zeros for a given solution depend on
initial conditions or boundary conditions for this solution and
are, in any case, easy to locate and compute.

There are several ways to ensure that a solution $y_n(x)$ has at most one zero in an interval; among them:

\begin{theorem} If one of the following conditions is satisfied in an interval $I$ (where all the coefficients
of the DDEs are continuous)
then $y_n$ and ($y_{n_i}$) have at most one
zero in the interval $I$ (trivial solutions excluded):
\begin{enumerate}
\item{$d_{n_i}(x) e_{n_i}(x)\ge 0$ in} $I$ \cite{Seg2}.
\item{$|\eta_i (x)|\ge 1$} in $I$ \cite{Seg3}.
\item{$\tilde{\mbox{A}}_n<0$} ($\tilde{\mbox{A}}_{n_i}<0$) \cite{Seg3}.
\end{enumerate}
\label{osccond}
\end{theorem}

The condition $e_{n_i}d_{n_i}<0$ is required for the method to apply. Furthermore, it is known that
when the DDEs (\ref{ddr}) are general, $d_{n_i} e_{n_i}$ can not change sign. Therefore $d_{n_i} e_{n_i}<0$
is a clear signature for the oscillatory character of the differential equation.

\subsection{Hypergeometric functions; selection of the optimal DDEs}

For hypergeometric functions several DDEs are available
for the construction of fixed point iterations (FPIs), depending on the
selection of contrast function.

Let us start by considering, for example, the case of the confluent
hypergeometric
 equation
\begin{equation}
xy^{\prime\prime}+(c-x)y^{\prime}-a y=0
\label{11}
\end{equation}

One of the solutions of this differential equation
are Kummer's series
$$M(a,c,x)\equiv  _1\!\!\mbox{F}_1(a;c;x)=\sum_{k=0}^{\infty}\Frac{(a)_k}{(c)_k k!}x ^k\,,
$$
for which different difference-differential relations are
available. Indeed, denoting $\alpha_n=\alpha+k\,n$,
$\gamma_n=\gamma+m\, n$ and $y_n\equiv M(\alpha_n ,\gamma_n ,x)$
we will have different sets of DDES (Eq.\ (\ref{ddr})) for
different selections of $(k,m)$.

For Gauss hypergeometric functions $_2\mbox{F}_1 (a,b;c;x)$, which are solutions
of the ODE
\begin{equation}
x(1-x)y''+[c-(a+b+1)x]y'-aby=0\,,
\label{GHF}
\end{equation}
\noindent the possible DDEs are determined by three-vectors with
integer components, that is, we will consider $y_n\equiv \:
_2\mbox{F}_1 (\alpha+k\,n,\beta+l\,n; \gamma+m\,n;x)$ and the
associated DDEs will be named $(k,l,m)$-DDEs. Finally, for the
case of the hypergeometric functions $_0\mbox{F}_1(;c;x)$, we can
only consider families $y_n= \:_0\mbox{F}_1 (;\gamma+k\,n ;x)$ and
the different relations are described by the integer numbers $k$.

Our FPMs can only be applied to solutions of second order ODEs.
This restricts our study to the hypergeometric functions
$_0\mbox{F}_1(;c;x)$,
 $_2\mbox{F}_0(a,b;;x)$,
 $_1\mbox{F}_1 (a;c;x)$ and
$_2 F_1 (a,b;c;x)$.

Regarding the selection of the different DDEs available, we will restrict ourselves to:
\begin{enumerate}
\item{DDEs} with continuous coefficients except at the singular points of the defining differential
equations.
\item{The} most simple DDEs in a given recurrence direction which allow the use of improved iteration steps.
Taking as example the case of confluent hypergeometric
functions, this means that the $(1,0)$-DDE will be described, and the
analysis of the $(-1,0)$, $(2,0)$,... DDEs will be skipped.
\end{enumerate}

The first restriction is convenient for simplicity and it means that the problem function and the contrast function
have zeros interlaced in each subdivision of the real interval defined by the singular points of
the differential equation; this is a convenient property for a simple application
of the FPMs and enables the application of each DDE to compute all the zeros in the different subintervals
of continuity of the solutions of the differential equation.

Regarding the second restriction, and considering the case of
confluent hypergeometric functions as example, it should be noted
that for the $(k,m)$-DDE,
generally two FPIs are available, one of them based on the ratio
$H_{-1}=y_n(x)/ y_{n-1}(x)$ and a second one based on
$H_{+1}(x)=y_n (x)/y_{n+1}(x)$. As described in \cite{Seg3},
generally one of these two iterations is preferable because
improved iteration steps can be considered according to Theorem
\ref{T2.4} (Theorem 2.4 in \cite{Seg3}). If we considered the
$(-k,-m)$-DDE, the two associated ratios $H_i$, $i=\pm 1$, would
be the same as before (replacing $i$ by $-i$). Because both
selections of DDEs are equivalent, only one of them will be
discussed.
By convention, we will consider pairs $(k,m)$ for which
the iteration  on $H_{-1}$ can take advantage of the monotony
property of $\hat{\mbox{A}}_n (x(z))$, as described in
\cite{Seg3}, for classical orthogonal polynomial cases (Jacobi, Hermite,
Laguerre).
 Once we have fixed this criterion the index $i$ in
Equations (\ref{ratH})-(\ref{ODEsnor}) can be dropped.
We consider the following additional notation: given a vector
$\vec{u}$ with integer components, we will denote by
DDE$(\vec{u})$ (FP$(\vec{u})$) the corresponding DDE (FPI) based
on the ratio $y_n/y_{n-1}$ for $x>0$.

On the other hand, and considering the confluent case as illustration, we will not analyze DDE$(2,2)$ nor any
successive multiples of the DDE$(1,1)$. This is so because the first restriction is generally
violated if successive multiples of a DDE are considered (there are exceptions to this; see the case of $_0\mbox{F}_1$
hypergeometric functions).

\subsubsection{Selection of the optimal DDEs}

There are several (and related) criteria to select among the available DDEs to compute the zeros
of a given function $y$. The associated FPMs tend to be more efficient as we are
closer to any of the following two situations:

\begin{enumerate}
\item{$\eta (x)=0$,}
\item{the} coefficient $\tilde{\mbox{A}}_n$ is constant.
\end{enumerate}

Of course, the first condition implies the second one (see Eq.\
(\ref{ODEsnor})). The first condition makes the FPI converge with
one iteration for any starting value. The second condition makes
the method an exact one using improved iteration steps (there is
even no need to iterate the FPM).

Let us recall that the FPIs associated to a given system of DDEs
are quadratically convergent to a zero $z_0$ (in the transformed
variable $z$) with asymptotic error constant $\eta(x(z_0))$.
Therefore, the smaller $|\eta (x)|$ is the fastest the convergence
is expected to be, at least for starting values close enough to
$z_0$.  On the other hand, the smaller the absolute value of
variation of $\tilde{\mbox{A}}_n$ is, the better the improved
iteration (Theorem \ref{T2.4}) will work because this implies the
exactness of the iteration criteria to estimate starting values
from previously computed zeros. This second criterion (on the
variation of $\tilde{\mbox{A}}$) is more difficult to apply, as we
will later see. It is, however, more relevant to improve the
iterative steps for obtaining starting values to compute zeros
than to improve the local convergence properties of the fixed
point methods, which are quadratically convergent anyway.

Indeed, as was described in \cite{Gil}, the natural FPIs for
orthogonal polynomials of confluent hypergeometric type
(FP$(-1,0)$) tend to converge slowly for the computation of the
first positive zeros when they become very small. This, for
instance, is the case  for Laguerre polynomials $L_{n}^{\alpha}
(x)$ when $\alpha\rightarrow -1^+$. The reason for this behavior
lies in the fact that the associated change of variables is
singular at $x=0$:
$$
z=\displaystyle\sqrt{(b-a)(1-a)}\ln x\,,
$$
In this way, the interval of orthogonality
 for the Laguerre polynomials $(0,+\infty)$
is transformed into $(-\infty,\infty)$ in the $z$ variable. This
means that the zeros which are very small in the $x$ variable,
tend to go to $-\infty$ in the $z$ variable. Therefore after
computing the second smallest zero, $x_2$, the next initial guess
for the FPI, $z(x_2)-\pi /2$, may lie well far apart for the value
$z(x_1)$, $x_1$ being the smallest zero.
 Although it is guaranteed
that the FPI will converge to $z(x_1)$, it could take a considerable
number of iterations to approach this value.

Fortunately, we will see that the rest of FPMs (different to
FP$(k,0)$) do not show such a singularity; therefore, we expect
better behavior near $x=0$ for these iterations.

This suggests that, given two FPIs with associated change of
variables $z_1(x)$, $z_2 (x)$ respectively, one should choose that
one which gives the largest displacement in the $x$ variable for
the same step in the corresponding $z$ (the typical value being
$\pi /2$). Let us stress that the possibility of passing the next
zero is ruled out: in the algorithms the sequence of all $z$
values calculated in a backward (forward) sweep form monotonically
decreasing (increasing) sequences.

We will therefore say that the change of variables $z_1$ behaves better than $z_2$
if $x(z_1+\Delta z)>x(z_2+\Delta z)$, for a typical value of $\Delta z$ ($\approx \pi /2$). Given the
definition of the changes of variables $z(x)$:
$$
x(z+\Delta z)-x (z)=\int_{z}^{z+\Delta z}\Frac{1}{\tilde{d}_{n}(x(z))}dz
$$
\noindent where $\tilde{d}_n=\sqrt{-d_n e_n}$, we can say that the
change of variable $z_1 (x)$ (and its associated FPI) is more
appropriate than the change $z_2 (x)$ when its coefficient
$\tilde{d}_{n}=\sqrt{-d_n e_n}$ is smaller than the corresponding
coefficient for $z_2(x)$.

Therefore, an alternative non-local prescription to that one
dealing with $\tilde{\mbox{A}}_n$ is the following: among the
possible DDEs and associated fixed point iterations, choose that
one for which $|d_n e_n|$ is smallest. As we will see, this is an
easy to apply criterion which correctly predicts the more
appropriate DDEs and FPI depending on the range of the parameters
and the dependent variable.

\section{Analysis of Hypergeometric functions}

We will use the DDEs satisfied by hypergeometric series as generated by
the Maple package {\it hsum.mpl} \cite{Koe}.
The results for each of the family of functions (change of
variable $z(x)$, function $\eta(x)$, etc.) considered can be
automatically generated using the package {\it zeros.mpl}
\cite{Gil}.

\subsection{Hypergeometric function $_0\mbox{F}_1 (;c;x)$}

The ODE satisfied by the function $y(x)= \:_0\mbox{F}_1(;c;x)$ is
\begin{equation}
x^2 y'' +c x y' -xy=0
\label{01}
\end{equation}
The solutions of these differential equations have an infinite
number of zeros for negative $x$ and are related to Bessel
functions:
$$
_0\mbox{F}_1 (;c;z)=\Gamma (c) (-z)^{(1-c)/2}J_{c-1}(2\sqrt{-z})
$$

\subsubsection{First DDE}

Let us consider the DDEs for the family of functions $y_n = \:
_0\mbox{F}_1(;\gamma +n;-x)$

The DDEs for this family
reads:
\begin{equation}
\begin{array}{l}
y'_{n}=-\Frac{c-1}{x}y_{n-1}+\Frac{c-1}{x}y_n\\[2mm]
y'_{n-1}=-\Frac{1}{c-1}y_{n} \label{DDE01}
\end{array}
\end{equation}
\noindent
where $c=\gamma+n$.

The relation with Bessel functions can be expressed saying that, if $y (x)$ is a solution of
 $x^2 y'' +(\nu+1) x y' +xy=0$ then
\begin{equation}
w(x)=x^{\nu} y(x^2 /4)
\label{relac}
\end{equation}
\noindent is a solution of the Bessel equation $x^2 w'' +x w' +(x^2-\nu^2 )w=0$.

Transforming the DDEs (\ref{DDE01}) as described in \cite{Seg3},
and in Section \ref{repaso}, the relevant functions are:
\begin{equation}
\eta (z(x)) =\Frac{\nu -1/2}{2\sqrt{x}},\\
\tilde{\mbox{A}}_n=1-\Frac{\nu^2 -1/4}{4 x},\\
z(x)=2\sqrt{x} \;
\end{equation}

The fixed point method deriving from this set of DDEs will have
identical performance as the system considered in \cite{Seg3},
which holds for Ricatti-Bessel functions
$j_{\nu}(x)=\sqrt{x}J_{\nu} (x)$, that are solutions of the
 second
order ODEs $y''+A(x)y=0$, with $A(x)=1-(\nu^2-1/4)/x^2$. The
identification of both methods with the replacement $x\rightarrow
x^2/4$ (according to Eq.\ (\ref{relac}) and to the change of
variable $z(x)$) is evident by comparing the $A(x)$
 and $\tilde{\mbox{A}}_n (x)$ coefficients. This is not surprising given that both methods compare the same problem
function, $J_{\nu}(x)$, with the same contrast function, $J_{\nu -1}(x)$, up to factors which do not vanish
(for example, the factor $\sqrt{x}$ for Ricatti-Bessel functions) and up to changes of variable.

\subsubsection{Second DDE}

For this type of hypergeometric functions, the only alternative
DDEs that can be built are those based on the family of functions
$y_n=\:_0\mbox{F}_1(;\gamma+m\,n;-x)$; $m=1$ corresponds to the
DDEs (\ref{DDE01}). We will only consider the case $m=2$
(equivalent, in the sense described before, to $m=-2$). For
$|m|>2$, the DDEs violate the first imposed condition on the
continuity of the coefficients. For the functions
$y_n=\:_0\mbox{F}_1 (;\gamma+2n;-x)$, the associated DDEs read:
\begin{equation}
\begin{array}{l}
y'_{n}=-\Frac{(c-2)^2+(c-2)-x}{(c-2) x}y_{n}+\Frac{c-1}{x}y_{n-1}\\
\\
y'_{n-1}=-\Frac{1}{c-2}y_{n-1}-\Frac{x}{(c-1)(c-2)^2}y_{n}
\label{DDE012}
\end{array}
\end{equation}
\noindent
where $c=\gamma+2n$.
 The relevant functions in this case are (again, writing $\nu=c-1$):
\begin{equation}
\begin{array}{l}
\eta(z(x))=\Frac{(\nu-1)^2}{2x}-1\\
\\
\tilde{\mbox{A}}_n(z(x))=\left(\Frac{\nu-1}{2x}\right)^2 (4x-(\nu^2-1))\\
\\
z(x)=\Frac{x}{\nu -1}
\end{array}
\end{equation}

\subsubsection{Comparison between DDEs}

The second DDE is no longer equivalent to the first one; in fact,
it has quite different characteristics. Let us compare the
expected performance of these two DDEs, according to the different
criteria described above.

To begin with, the $\eta(x)$ parameter never vanishes for the
first DDE (DDE1), except when $\nu=1/2$, in which case
 the method
with improved steps is exact without the need to iterate the FPI
even once (forward or backward sweeps \cite{Seg2} are used
depending on the sign of $\eta (x)$). In contrast, DDE2 has an
$\eta$-function which changes sign at $x_{\eta}=(\nu -1)^2/2$ and
the zeros are computed by an expansive sweep \cite{Seg2}. Close to
$x_{\eta}$ we can expect that DDE2 tends to behave better in
relation to local convergence, because the asymptotic error
constant tends to be small.

We observe that, as $x\rightarrow +\infty$,  $\eta(x)$ goes to zero for DDE1 but it tends to $-1$ for DDE2. This suggests
that DDE1 will have faster local convergence than DDE2 for large $x$. On the other hand, as $\nu$ increases, $\eta(x)$
becomes larger, however, it is difficult to quantify the impact on local converge because as $\nu$ increases also
the smallest zero becomes larger. Let us also take into account that DDE2 will have small $\eta (x)$ for
$x$ close to $x_{\eta}$, which becomes large for large $\nu$.

Regarding the behavior of $\tilde{\mbox{A}}_n (z)$, it is
monotonic for DDE1 and has a maximum for DDE2, which allows the
use of improved iteration steps. For DDE1, $\tilde{\mbox{A}}_n
(z)$  is constant when $\nu=1/2$, which means that sweep with
improved iteration steps is exact, as commented  before. The
maximum for DDE2 is at $x_m=(\nu^2-1)/2$, where
$\tilde{\mbox{A}}_n(x_m)=(\nu-1)/(\nu+1)$. Around this extremum,
the improved iteration steps tend to work better because
$\tilde{\mbox{A}}_n (x)$ will be approximately constant; how
constant $\tilde{\mbox{A}}_n (x)$ is around $x_m$ can be measured
by the convexity at this point. We find:
$$
\ddot{\tilde{\mbox{A}}}_n (z_m)=-8\Frac{\nu -1}{(\nu +1)^3}
$$
\noindent where $z_m=z(x_m)$. As $\nu$ becomes larger, $\ddot{\tilde{\mbox{A}}}_n (z)$  
becomes smaller around the maximum of $\tilde{\mbox{A}}_n (z)$ and 
the improved iteration will work better. This fact again, favors DDE2 for large $\nu$.

Finally, considering the criterion of smaller $D_n=|d_n e_n|$, we find that, for DDE1
\begin{equation}
D_n=\Frac{1}{x}
\end{equation}
\noindent while for DDE2
\begin{equation}
D_n=\Frac{1}{(\nu -1)^2}
\end{equation}

This again shows that DDE1 will improve as $x$ increases while DDE2 will be better for large $\nu$. Numerical
experiments show that for $\nu>100$ the second DDE is preferable over the first, particularly for computing
the smallest zeros.

The different criteria yield basically the same information. However the prescription on $D_n=|d_n e_n|$ is the simplest
one to apply. From now on, we will not repeat the analysis for the different criteria. Instead, we adopt this last
criterion to analyze the rest of cases.

\vspace*{0.5cm}
\begin{minipage}{6cm}
\centerline{\protect\hbox{\psfig{file=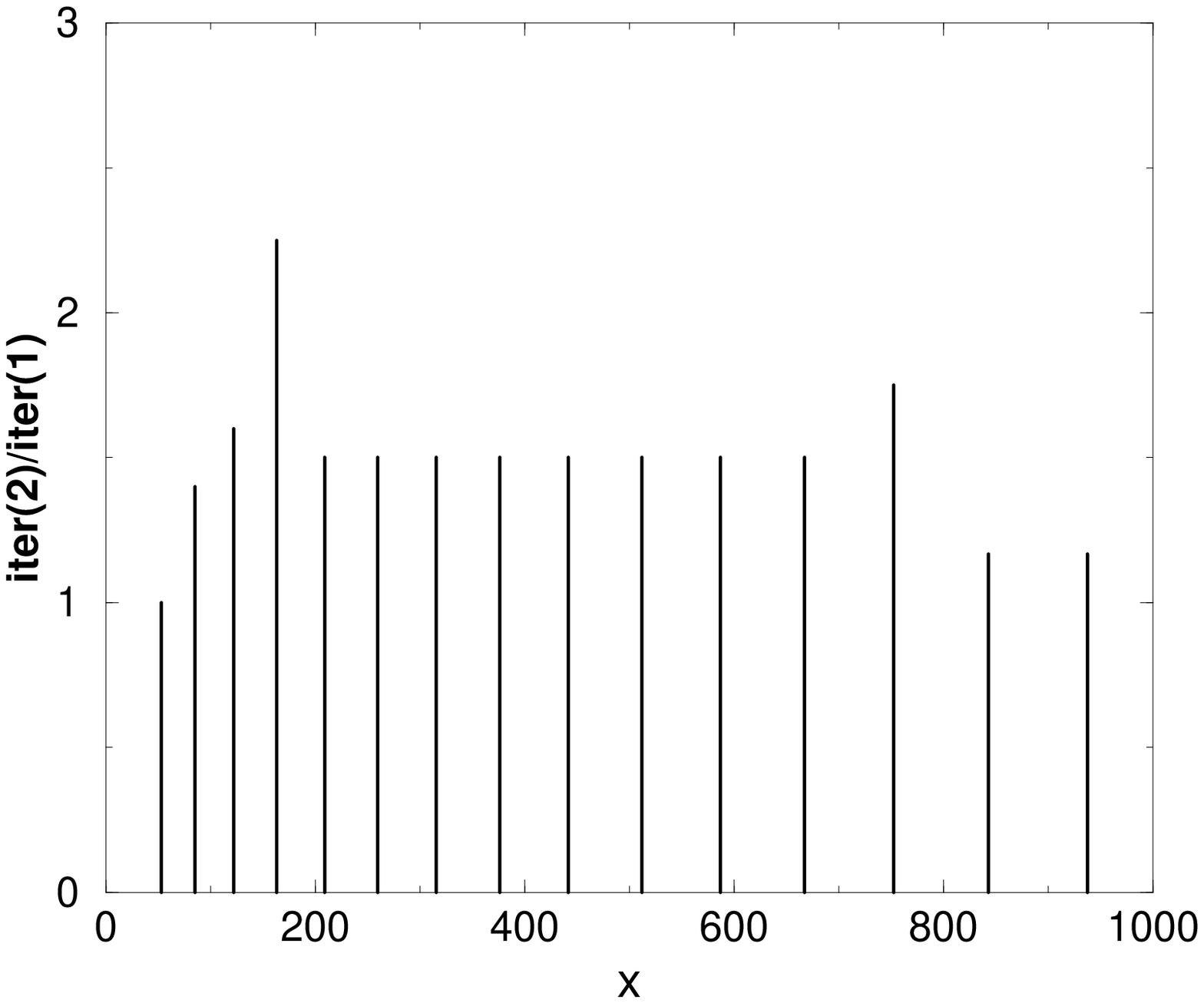,width=5.cm}}}
\end{minipage}
\hfill
\begin{minipage}{6cm}
\centerline{\protect\hbox{\psfig{file=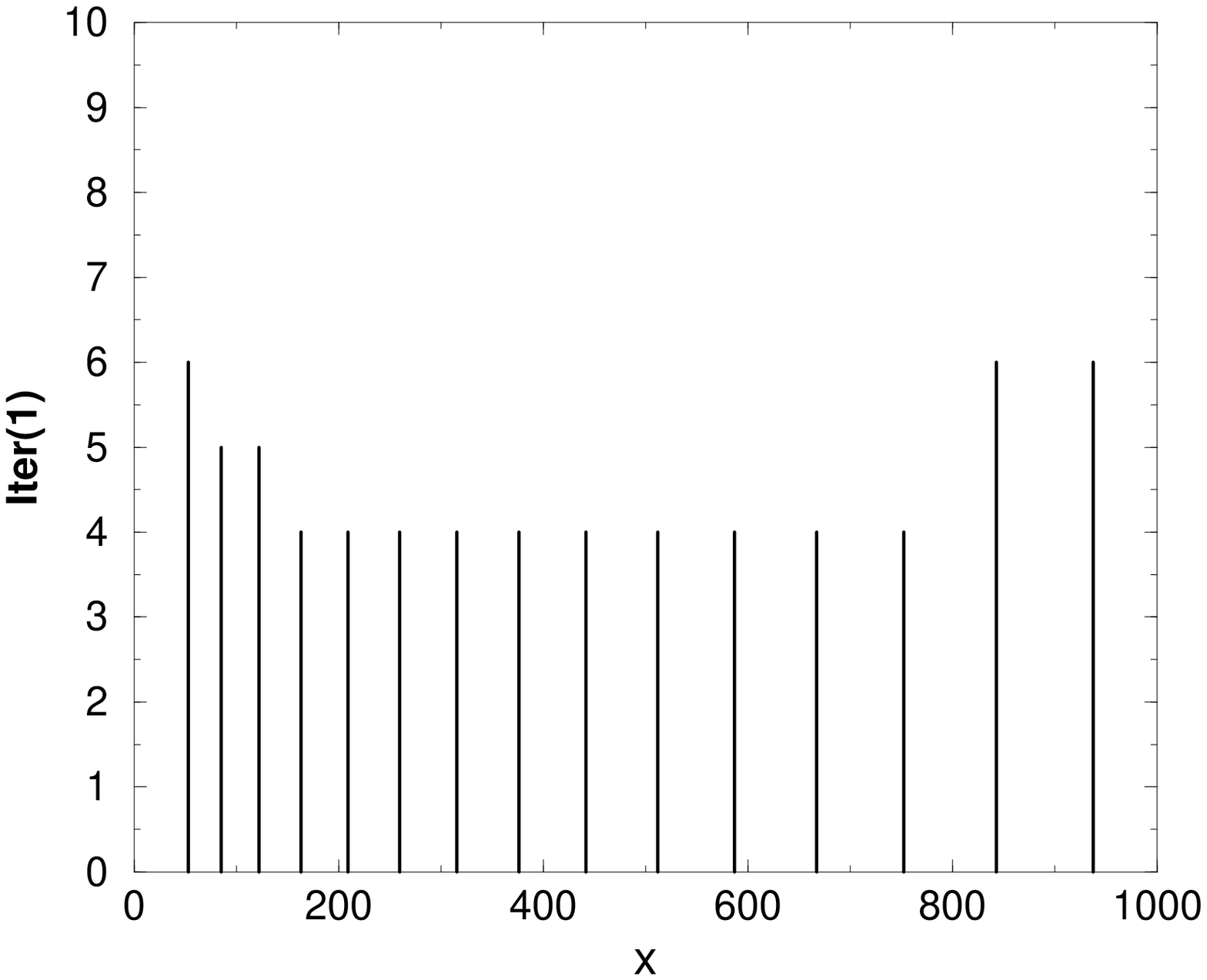,width=5.cm}}}
\end{minipage}
{ {\bf Figure 1.} {\bf Left}: Ratio between the number of
iterations needed for the second and first DDEs for the
computation of the zeros of $_0\mbox{F}_1(;11,-x)$ (the zeros of
$J_{10}(2\sqrt{x})$). {\bf Right}: Number of iterations needed for
the first DDE. }

\begin{minipage}{6cm}
\centerline{\protect\hbox{\psfig{file=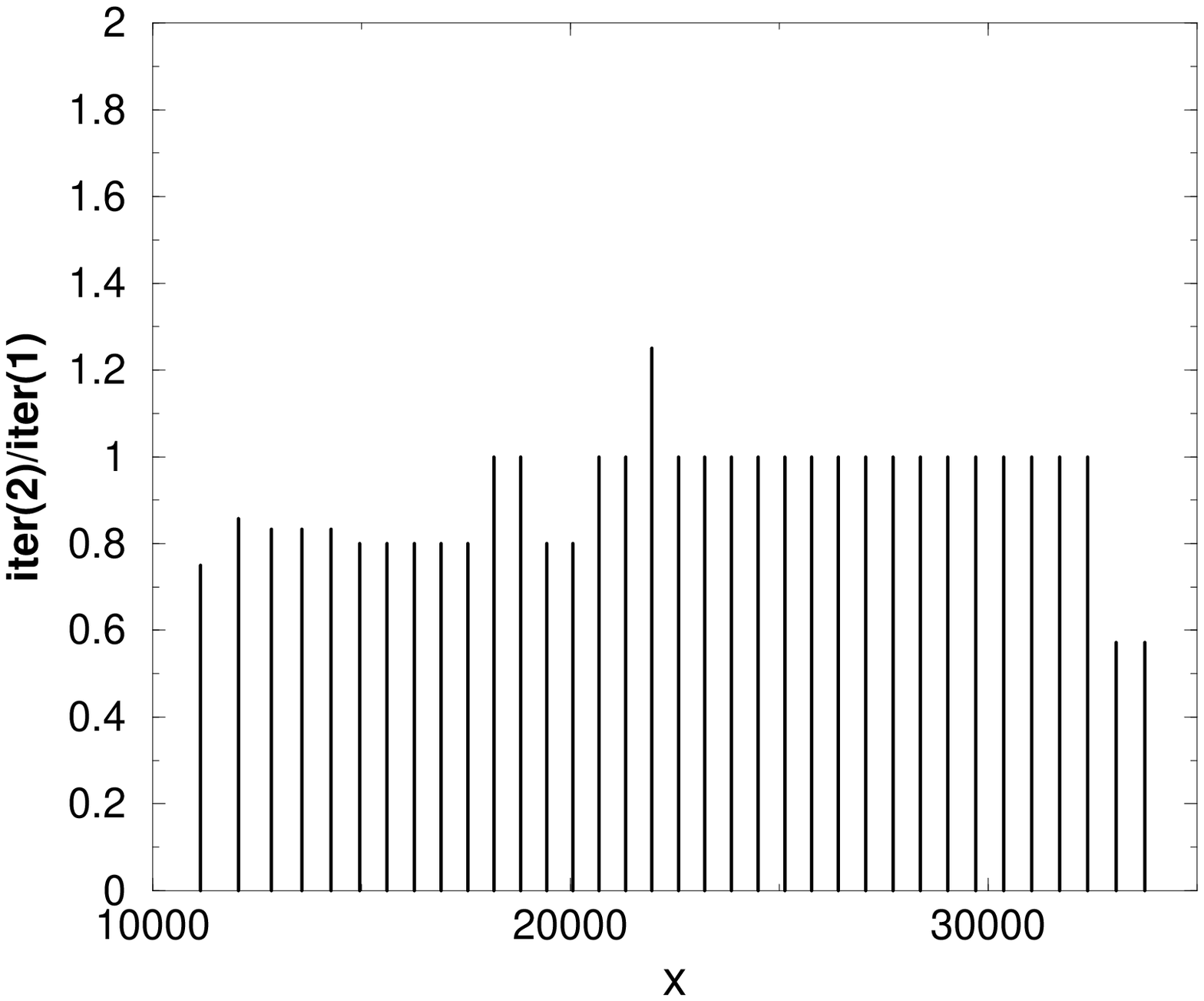,width=5.cm}}}
\end{minipage}
\hfill
\begin{minipage}{6cm}
\centerline{\protect\hbox{\psfig{file=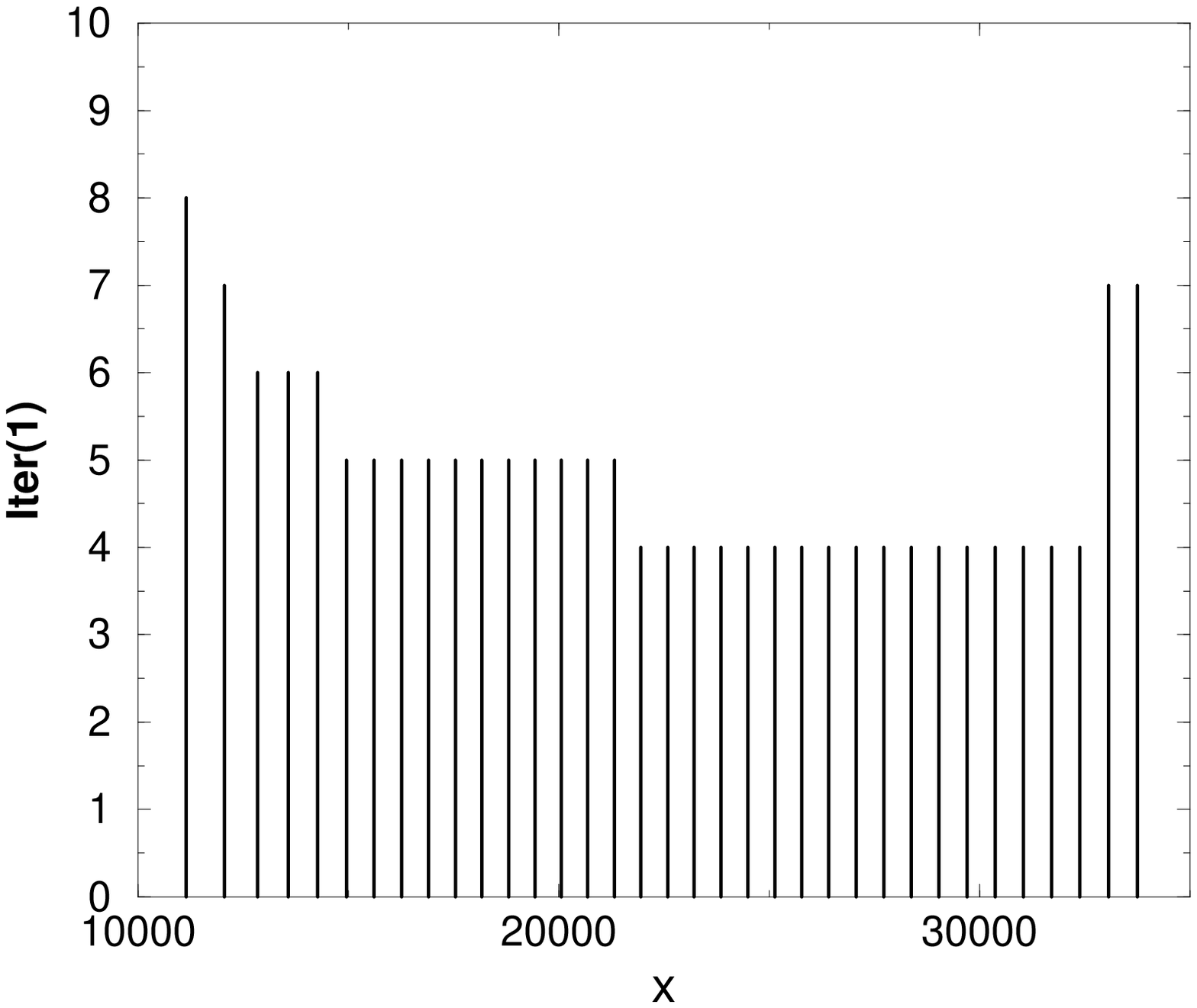,width=5.cm}}}
\end{minipage}
{ {\bf Figure 2.}
{\bf Left}: Ratio between the number of iterations
needed for the second and first DDEs for the computation of the zeros
of $_0\mbox{F}_1(;201,-x)$ (the zeros of $J_{200}(2\sqrt{x})$).
{\bf Right}: Number of iterations needed for the first DDE.
}

\subsection{Confluent hypergeometric function}

For the confluent hypergeometric case we have a larger variety of
DDEs to choose, because we can choose families of functions
$y_n=\:_1\mbox{F}_1(a+k\, n; c+m\, n;x)$, or, more generally,
$y_n=\phi (a+k\, n; c+m\, n;x)$, being $\phi$ any solution of Eq.\
(\ref{11}). The families which give rise to DDEs satisfying all
our requirements are three, corresponding to the following
selections of $(k,m)$: $(1,0)$, $(1,1)$ and $(0,-1)$.

We will give the corresponding DDEs and the associated functions. At the same time, we will
restrict the range of parameters for which the functions are oscillatory (considering the
first oscillatory condition in Theorem \ref{osccond}).

We can restrict the study to $x>0$ because, if $\phi (a;c;x)$ is a
solution of Eq.\ (\ref{11}), then $e^{x}\phi (c-a;c;-x)$ is also a
solution of Eq.\ (\ref{11}).

 \subsubsection{$(k,m)=(1,0) \rightarrow y_n=\phi$($\alpha +n$; $\gamma$ ; x)}

Let us write, for shortness and in order to
compare with other recurrences $a=\alpha+n$, $c=\gamma$. As before commented,
we consider simultaneously the equivalent directions $(k,m)=(1,0)$ and $(k,m)=(-1,0)$ but
we present only the DDEs and related functions for the recurrence direction for
which the FPI based on the ratio $H_{-1}=y_n /y_{n-1}$ can be used with improved
iteration steps (Theorem \ref{T2.4}); this is the direction $(1,0)$.

The $(-1,0)$ direction is the natural one for orthogonal polynomials of hypergeometric type
(Laguerre, Hermite), which
are related to confluent hypergeometric series of the
type $_1\mbox{F}_1(-n;\gamma ;x)$ (see, for instance, \cite{Leb}, Eqs. (9.13.8-10)).

The DDEs for $(k,m)=(1,0)$ read
\begin{equation}
\begin{array}{ll}
y'_n=&\Frac{a-c+x}{x}y_n-\Frac{a-c}{x}  y_{n-1}\\
\\
y'_{n-1}= &-\Frac{a-1}{x}y_{n-1}+\Frac{a-1}{x}y_n
\label{10}
\end{array}
\end{equation}
\noindent and therefore, applying the first oscillatory condition of Theorem \ref{osccond},
 the parameters are restricted to:
$$
(a-c)(a-1) > 0
$$
\noindent otherwise the functions $y_{n}$, $y_{n-1}$ will be non-oscillatory. If we repeat the same
argument for $(k,m)=(-1,0)$ or, equivalently, apply the same criteria for the DDEs for the functions
$y_n$, $y_{n+1}$ the following restriction is obtained:
$$
(a-c+1)a >0\,.
$$

The associated functions for this DDE are
\begin{equation}
\begin{array}{l}
\eta (z(x))=-\Frac{2a+c+x}{2\displaystyle\sqrt{(c-a)(1-a)}}\\
\\
\tilde{\mbox{A}}_n(z(x))= \Frac{-x^2+2(c-2a)x-(c-1)^2}{4(c-a)(1-a)}\\
\\
z(x)=\displaystyle\sqrt{(c-a)(1-a)}\ln x
\end{array}
\end{equation}

Let us notice that $\eta (x)$ becomes negative for large $x$,
which gives the most appropriate sweep for large $x$ (forward)
since $\tilde{\mbox{A}}_n (z(x))$ decreases for large $x$.

Observe that the lack of a singularity in  $\tilde{\mbox{A}}_n (z(x))$ is only apparent because the function $z(x)$ is singular
at $x=0$.

\subsubsection{$(k,m)=(0,-1)\rightarrow y_n=\phi$($\alpha$; $\gamma -n$ ; x)}

As before, we denote $a=\alpha$, $c=\gamma -n$.  The DDEs for the system read:
\begin{equation}
\begin{array}{ll}
y'_n=& y_n  + \Frac{a-c}{c} y_{n-1}\\
\\
y'_{n-1}= & -\Frac{c}{x}  y_{n-1}  +\Frac{c}{x}y_n
\label{0-1}
\end{array}
\end{equation}
\noindent which implies the restriction $(c-a)x>0$ for the
solutions to have oscillatory nature. Repeating the same for
$(k,m)=(0,1)$, we arrive at the condition $(c-a-1)x>0$, which for
$x>0$ gives a more restrictive condition $c-a>1$ (for $x<0$, the
first condition is more restrictive and gives $a-c<0$).

Let us recall that the condition $c-a>1$ for $x>0$ means that, if this condition is not met,
neither $y_n$ nor $y_{n-1}$ can have two zeros in $x>0$. Given that we are interested in computing
zeros of oscillatory functions, we consider these restrictions. Let us however notice that the
possible isolated zero of $y_n$ for $x>0$ and $0<c-a<1$, could be also computed by means
of the FPM associated to the DDEs (\ref{0-1}).

Considering also the restrictions imposed in the previous selection of
DDEs, we obtain the following:
\begin{theorem}
Let $y$ be a solution of the confluent hypergeometric equation (\ref{11}) for $x>0$. If $y$
has at least two zeros then $c-a>1$ and $a<0$.

Let $y$ be a solution of the confluent hypergeometric equation (\ref{11}) for $x<0$. If $y$
has at least two zeros then $c-a<0$ and $a>1$.
\label{condiex}
\end{theorem}
For more detailed results on the number of zeros of confluent
hypergeometric functions, we refer the reader to \cite{Bat},
Volume 1, Section 6.16.

The associated functions are
\begin{equation}
\begin{array}{l}
\eta (z(x))=-\Frac{2c-1+2x}{4\displaystyle\sqrt{(c-a)x}}\\
\\
\tilde{\mbox{A}}_n(z(x))= \Frac{8c x-16x a-3+8c-4c^2-4x^2}{16(c-a)x}\\
\\
z=2\displaystyle\sqrt{(c-a) x}
\end{array}
\end{equation}

Let us notice that $\eta(z(x))$ becomes negative for large $x$,
which gives the most appropriate sweep for large $x$ (forward)
since, for positive $x$, $\tilde{\mbox{A}}_n (z(x))$ decreases for
large $x$.

\subsubsection{$(k,m)=(1,1)\rightarrow y_n=\phi$($\alpha +n$; $\gamma +n$ ; x)}

The DDEs are the following (writing $a=\alpha+n$, $c=\gamma +n$).
\begin{equation}
\begin{array}{ll}
y'_n=&   \Frac{x+1-c}{x}y_n +\Frac{c-1}{x}  y_{n-1}\\
\\
y'_{n-1}= & \Frac{a-1}{c-1} y_n
\end{array}
\label{11F11}
\end{equation}
\noindent which implies the oscillatory condition $(a-1)x>0$;
considering $(k,m)=(-1,-1)$ or, equivalently, the DDEs relating
$y_{n}$, $y_{n+1}$ and their derivatives, the condition obtained
is $ax>0$. These conditions are consistent with Theorem
\ref{condiex}.

The functions associated to these DDEs for $x>0$ are
\begin{equation}
\begin{array}{l}
\eta= -\Frac{2x+3-2c }{ 4 \displaystyle\sqrt{(1-a)x}}\\
\\
z=2\displaystyle\sqrt{(1-a) x}\\
\\
 \tilde{\mbox{A}}_n (z(x))= \Frac{16x a+4x^2-8x c+3-8 c+4 c^2}{16(-1+a)x}
\end{array}
\end{equation}

This iteration cannot be used for $c=1$ (see Eq.\ (\ref{11F11}))
unless the FPI stemming from the ratio $y_n/y_{n+1}$ is used, in
which case the improved iteration steps cannot be applied.

\subsubsection{Comparing fixed point iterations}

As commented before, we will consider the prescription consisting in choosing the DDEs for
which the product $D_n=-d_n e_n$ is smaller. Applying literally this criterion, we find the
following preferred regions of application:
\begin{enumerate}
\item{FP(1,1)} should be applied for $x<c-a$ and FP(1,0) for $x>c-a$.
\item{FP(1,1)} is always better than FP(0,-1).
\item{FP(0,-1)} is better than FP(1,0) for $x<1-a$, but $1-a<c-a$ and FP(1,1)
is better there.
\end{enumerate}

Therefore, the best combination of the considered FPIs is FP(1,1)
for $x<c-a$ and FP(1,0) for $x>c-a$. The iteration FP(0,-1) has
the same behavior as FP(1,1) and can be used as a replacement when
$c=1$ (in this case FP(1,1) can not be used). In fact, the
FP(0,-1) and FP(1,1) are not independent because, as it is well
known, if $\psi (\alpha ;\gamma ;x)$ are solutions of the
confluent hypergeometric equation
$x\psi''+(\gamma-x)\psi'-\alpha\psi =0$, then $y(a;c;x)=x^{1-c}
\psi (1+a+c;2-c;x)$ is a solution of $x y'' + (c-x)y' -a y=0$.

Let us illustrate this behavior with numerical examples:

In Fig.\ 3 we compare FP(1,1) with FP(1,0) for the case of
Laguerre polynomials $L_{n}^{\alpha}(x)$ for the quite extreme
case $n=50$, $\alpha=-0.9999$. On the left, the ratio between the
number of iterations employed by FP(1,1) and FP(1,0) is shown as a
function of the location of the zeros. The first two zeros are
skipped in the left figure to show in more detail the behavior for
large $x$ (for the first zero the ratio was 40 while for the
second it was 5). The improvement for small $x$ when considering
FP(1,1) is quite noticeable. For larger $x$ ($x>c-a\simeq 50$),
FP(1,0) works generally better than FP(1,1) but the improvement is
not so noticeable. On the right figure, the number of iterations
used to compute each zero when considering FP(1,1) is shown.

\begin{minipage}{6cm}
\centerline{\protect\hbox{\psfig{file=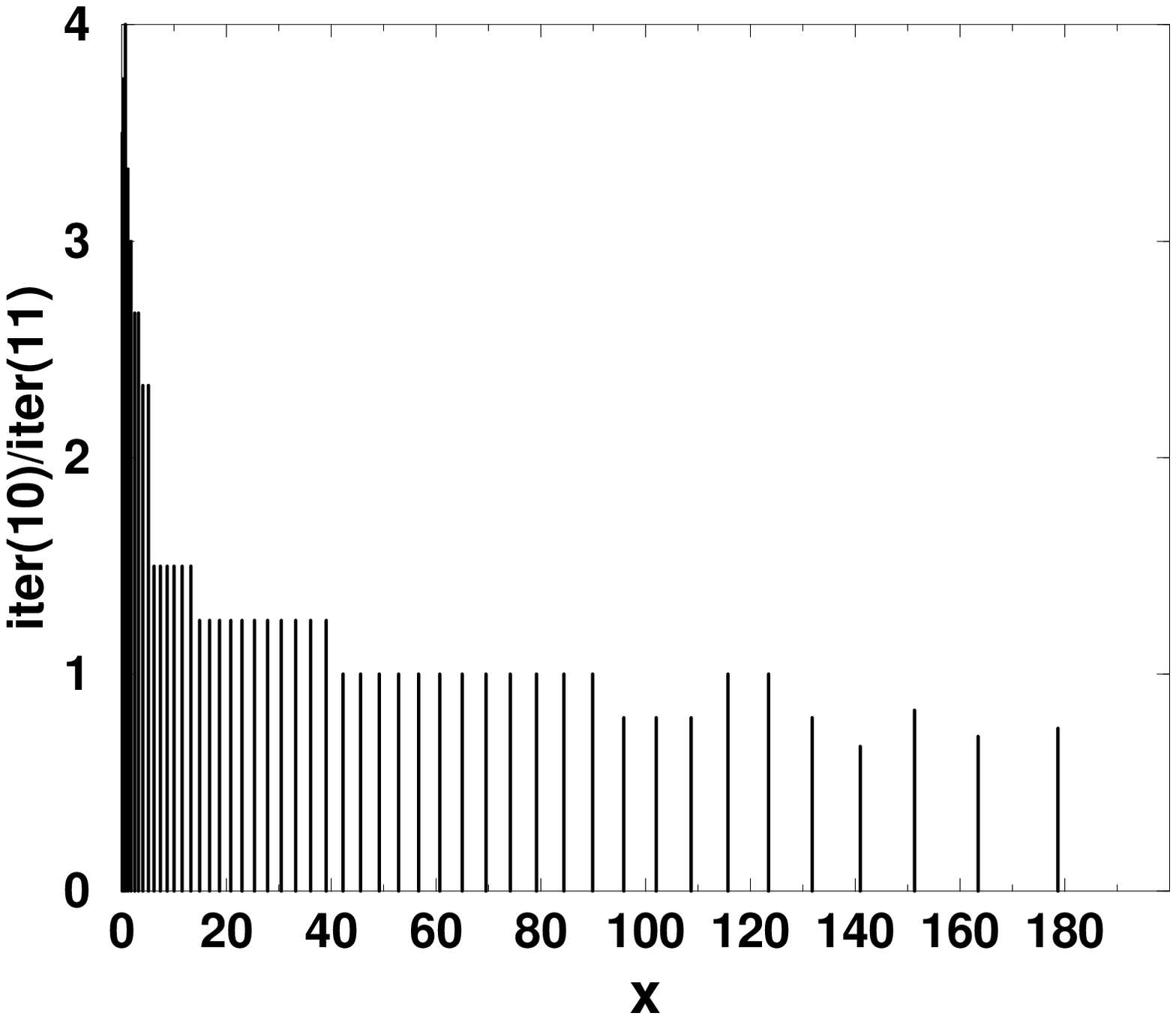,width=6.cm}}}
\end{minipage}
\hfill
\begin{minipage}{6cm}
\centerline{\protect\hbox{\psfig{file=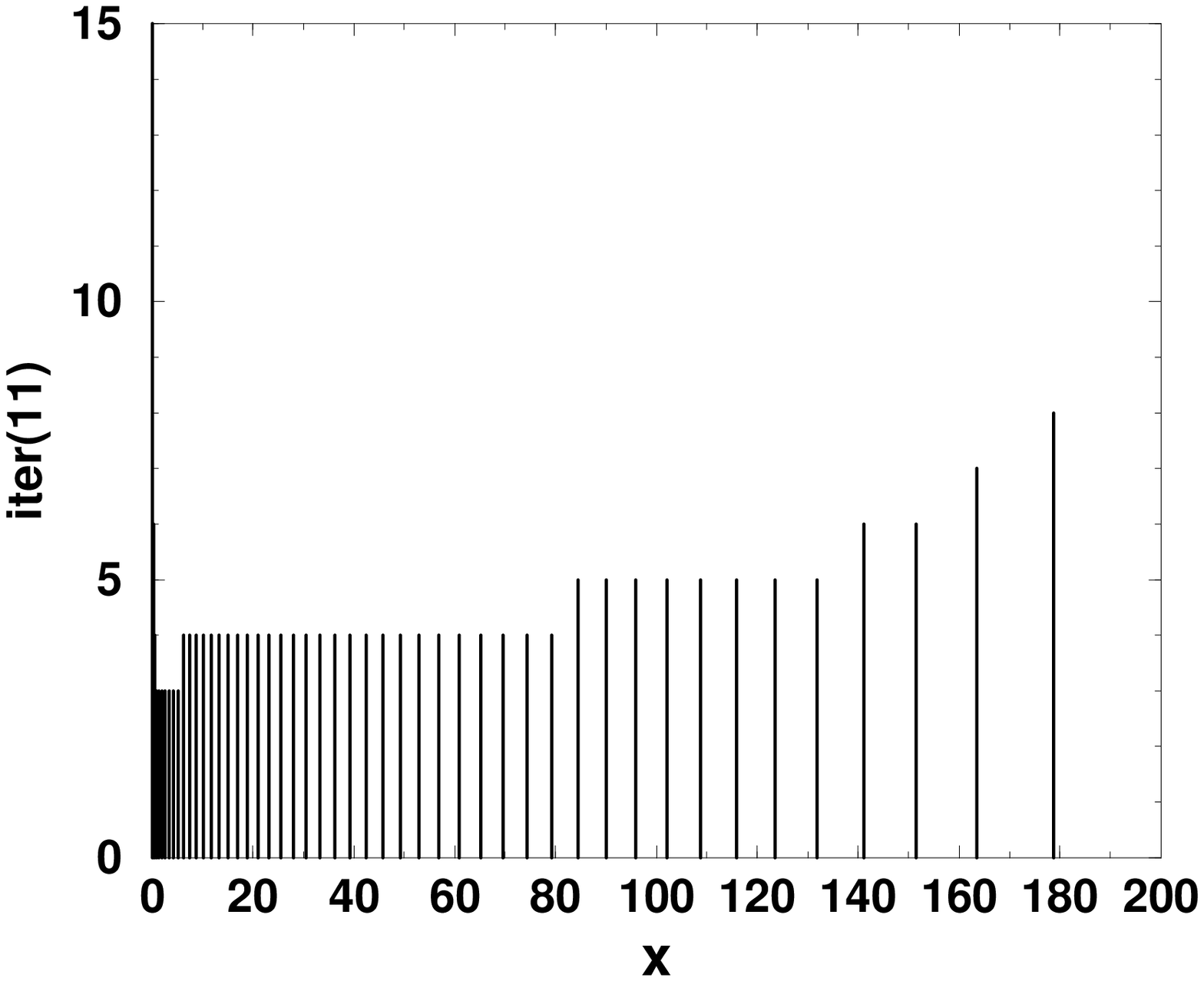,width=6.cm}}}
\end{minipage}

{ {\bf Figure 3.} {\bf Left:} Ratio between the number of iterations
needed by FP(1,0) and FP(1,1) for the calculation of the zeros of the generalized Laguerre polynomial
$L_{50}^{-.9999}(x)$, as a function of the location of the zeros.
{\bf Right:} Number of iterations needed by FP(1,1).}

In conclusion, FP(1,1) is a more appropriate choice than the natural recurrence for
orthogonal polynomials (FP(1,0)), although FP(1,0) slightly improves the convergence
of FP(1,1) for $x>c-a$.

In Fig.\ 4 we compare FP(1,0) with FP(0,-1) for the case of
generalized Laguerre polynomials but now for a choice of the
parameters $n=50$, $\alpha=0$. This situation corresponds to the
case $c=1$ where FP(1,1) can not be applied. As in Fig.\ 3, the
ratio of the number of iterations for the first two zeros is not
plotted (the ratio of the first zero was 8 and for the second it
was 5). As expected from our previous analysis, the iteration
FP(0,-1) behaves quite better than FP(1,0) for small $x$.

\begin{minipage}{6cm}
\centerline{\protect\hbox{\psfig{file=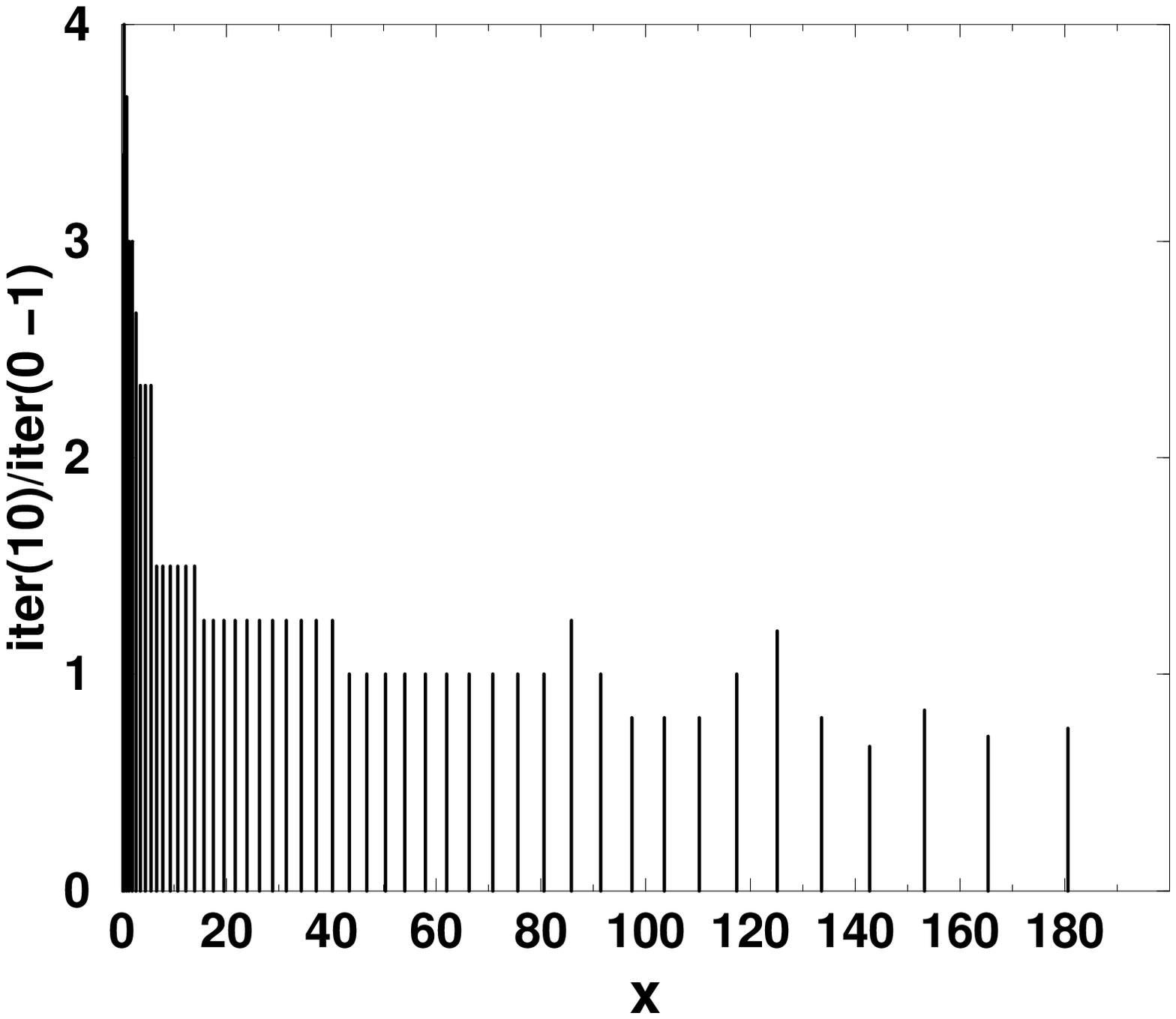,width=6.cm}}}
\end{minipage}
\hfill
\begin{minipage}{6cm}
\centerline{\protect\hbox{\psfig{file=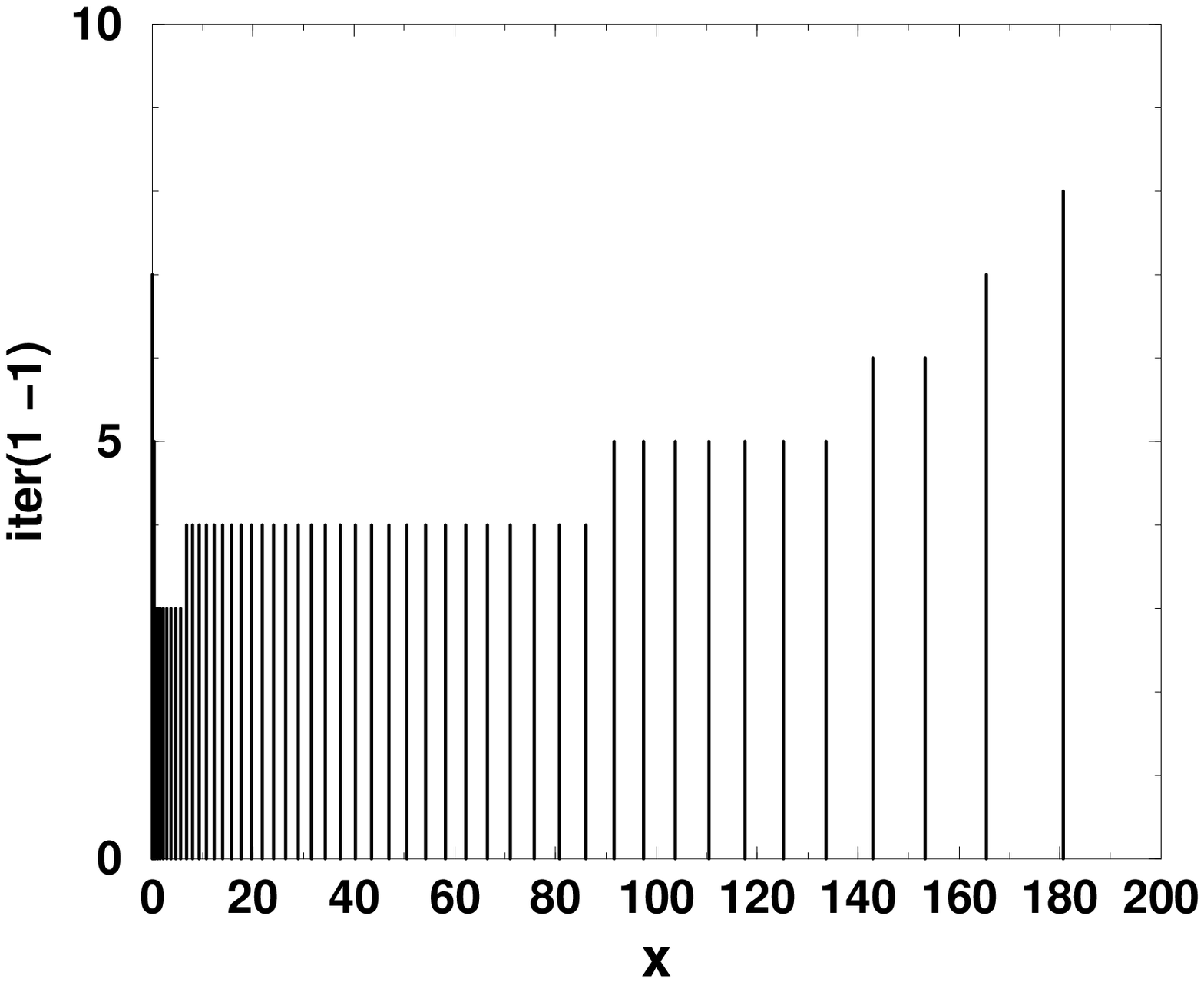,width=6.cm}}}
\end{minipage}

{{\bf Figure 4.} {\bf Left:} Ratio between the number of iterations
needed by FP(1,0) and FP(0,-1) for the calculation of the zeros of the generalized Laguerre polynomial
$L_{50}^{0}(x)$, as a function of the location of the zeros.
{\bf Right:} Number of iterations needed by FP(0,-1).}

\subsection{Mysterious hypergeometric function}

This is the name given to the hypergeometric series $_2\mbox{F}_0 (a,b;;x)$,
 which diverges for all $x\neq 0$ (except in the terminating cases)
and can only be interpreted
in an asymptotic sense. It is well know that, for negative $x$ we have:
\begin{equation}
_2\mbox{F}_0(a,b;;x) = (-1/x)^a U(a,1+a-b,-1/x)
\label{myst}
\end{equation}
\noindent where $U(a,c,x)$ is a solution of the confluent hypergeometric
equation (\ref{11}). The functions $_2\mbox{F}_0$ can
 be analytically continued
 to the whole complex
plane cut along the line $\Re(z)>1$, $\Im(z)=0$. The function (\ref{myst}) is
a solution
of the 2-0 hypergeometric differential equation:
\begin{equation}
x^2 y''+[-1+x(a+b+1)]y'+ab\,y=0\,.
\label{mysE}
\end{equation}

In general terms, without referring to any particular solution of
the corresponding differential equations, the problem of computing
the real zeros of the mysterious hypergeometric function for
negative $x$ can be transformed into a problem of computation of
the zeros of confluent hypergeometric functions Eq.\ (\ref{myst}).
This is so because one can check that, if we denote by
$y(\alpha,\beta,x)$ a set of solutions of the confluent
hypergeometric equation, then $w(x)=|x|^{-a} y(a,1+a-b,-1/x)$, for
$x>0$ or $x<0$, are solutions of Eq.\ (\ref{mysE}).

For this reason and for brevity we omit further details.

\subsection{Gauss Hypergeometric Functions}

   Let us consider the hypergeometric function $_2\mbox{F}_1(a,b;c;x)$. We will consider the
 DDEs for families of functions of the type $\psi (\alpha + kn, \beta +m n; \gamma + l n; x)$, with
 $\psi (a,b;c;x)$ solutions of the hypergeometric equation (\ref{GHF}); we use the DDEs for $_2\mbox{F}_1 (a,b;c;x)$
series as generated by ${\it hsum.mpl}$.

Similarly as we did for the confluent case, we can obtain oscillatory conditions
for the coefficients $a$, $b$ and $c$, depending on the range of $x$.
 If these conditions are not satisfied by the parameters,
then we can assure that if there exists one zero of the function, this is an isolated zero.  As before, these conditions
are obtained by requiring that $d_n e_n<0$; combining the restrictions imposed by this condition for the DDEs that
we will later show we obtain the following:

\begin{theorem} Let $\psi (a,b;c;x)$ be a solution of Eq.\ (\ref{GHF}) defined in $(-\infty, 0)$, then if this function
is oscillatory in this interval (it has at least two zeros) then one of the following sets of conditions must be verified:
  $$
   \left\{
   \begin{array}{l}
   \mbox{(C1):  } a<0\,,\,b<0\,,\,c-a>1\,,\,c-b>1 \\
    \mbox{(C2):  }  a>1\,,\,b>1\,,\,c-a<0\,,\,c-b<0
   \end{array}
   \right.
  $$

 Similarly,if $\psi (a,b;c;x)$ is a solution of Eq.\ (\ref{GHF}) defined in $(0, 1)$, the oscillatory conditions are
   $$
   \left\{
   \begin{array}{l}
    \mbox{(C3):  } a<0\,,\,b>1\,,\,c-a>1\,,\,c-b<0 \\
     \mbox{(C4):  }  a>1\,,\,b<0\,,\,c-a<0\,,\,c-b>1
   \end{array}
   \right.
  $$

 Finally, if $\psi (a,b;c;x)$ is a solution of Eq.\ (\ref{GHF}) defined in $(1, +\infty)$ and is oscillatory in this
interval,
 one of the following sets of conditions must be verified:
   $$
   \left\{
   \begin{array}{l}
    \mbox{(C5):  } a<0\,,\,b<0\,,\,c-a<0\,,\,c-b<0 \\
    \mbox{(C6):  } a>1\,,\,b>1\,,\,c-a>1\,,\,c-b>1
   \end{array}
   \right.
  $$
\label{oscG}
\end{theorem}

  The different sets of conditions for the three subintervals in Theorem \ref{oscG} can be obtained by combining
the restrictions obtained
from the following values of $(k,l,m)$: $(\pm 1,0,0)$, $(\pm 1,\pm 1,0)$, $(\pm 1, 0,\pm 1)$, $(\pm 1,\mp 1,0)$,
 $(0,0,\pm 1)$ and $(\pm 1,\pm 1, \pm 1)$. In any case, the analysis for the three different subintervals are
not independent because, as it is well know (see \cite{Bat}, Vol.\
I, Chap.\ II), if we denote by $\psi (\alpha,\beta;\gamma,x)$ the
solutions of the hypergeometric equation
$x(1-x)y''+(\gamma-(\alpha +\beta +1)x)y'-\alpha\beta y=0$ in the
interval $(0,1)$ one can write solutions in the other two
intervals by using that both
\begin{equation}
y(a,b;c;x)=(1-x)^{-a} \psi (a,c-b;c;x/(x-1))\,,x<0
\label{x<0}
\end{equation}
\noindent and
\begin{equation}
y(a,b;c;x)=x^{-a} \psi (a,a+1-c;a+b+1-c;1-1/x)\,,x>1
\label{x>1}
\end{equation}
are solutions of the hypergeometric differential equation
$x(1-x)y''+(c-(a+b+1)x)y'-aby=0$. With this, it is easy to see
that the conditions C1 and C2 can be obtained from C3 and C4
respectively by using Eq.\ (\ref{x<0}), while C5 and C6 derive
from C3, C4 and Eq.\ (\ref{x>1}).

 Notice, in addition, that the 6 conditions in Theorem \ref{oscG} are mutually exclusive, which means that

\begin{theorem}
Given three values of the parameters $a$, $b$ and $c$, at most one
of the subintervals $(-\infty,0)$, $(0,1)$, $(1,+\infty)$
possesses oscillatory solutions (solutions with at least two
zeros).
\end{theorem}

 The oscillatory conditions C6 is of no use for non-terminating hypergeometric series,
 because they diverge for $x>1$; however, this is a possible case for
 other solutions of the differential equation. The conditions C3 correspond
 to Jacobi polynomials $P_{n}^{(\alpha ,\beta)}(x)=\Frac{(\alpha +1)_n}{n!}$ $_
2\mbox{F}_1(-n,1+\alpha+\beta+n;\alpha +1;(1-x)/2)$ ($\alpha
,\beta >-1$) of order $n\ge 2$ (for order $n=1$ the conditions are
not satisfied, not surprisingly because according to our criteria
a polynomial of degree 1 is non-oscillating). Particular cases of
Jacobi polynomials are Gegenbauer ($\alpha=\beta$), Legendre
($\alpha=\beta=0$) and Chebyshev ($\alpha=\beta =-1/2$)
polynomials.

\subsubsection{DDEs and change of variables}

 In this section, we compile the expressions for the different DDEs as well as the associated
change of variable. For brevity, the associated functions $\eta
(x)$ and $\tilde{\mbox{A}}_n (x)$ are not shown. 

It is understood that $a=\alpha + k\,n$, $b=\beta +l\,n$, $c=\gamma +m\,n$ for DDE$(k,l,m)$.

\begin{enumerate}
\item{DDE$(1,0,0)$}
$$
\begin{array}{ll}
y'_n=&\Frac{-a-bx+c}{x(x-1)}y_n+\Frac{a-c}{x(x-1)}  y_{n-1}\\
\\
y'_{n-1}= & \Frac{1-a}{x}y_{n-1}+\Frac{a-1}{x}y_n\\
\end{array}
$$

Change of variable:
$z(x)=-2\sqrt{(c-a)(1-a)}\tanh^{-1}(\sqrt{1-x})$

\item{$DDE(1,1,0)$ :}
$$
\begin{array}{ll}
y'_n=&\left(\Frac{-(a+b)(a+b-c-1)+ab-c}{
x(a+b-c-1)}+\Frac{a+b-c}{x(1-x)}\right)y_n+
\Frac{(a-c)(c-b)}{x(1-x)(a+b-c-1)}
  y_{n-1}\\
\\
y'_{n-1}= & \Frac{a+b-ab-1}{(a+b-c-1)x} y_{n-1}
-\Frac{(1-x)(a-1)(1-b)}{(a+b-c-1)x}y_n\\
\end{array}
$$

Change of variable: $z(x)=\Frac{\sqrt{(b-c)(c-a)(b-1)(1-a)}}{|
a+b-c-1|}\ln x$,  $x>0$.

\item{$DDE(1,1,2)$ :}
$$
\begin{array}{ll}
y'_n=& \Frac{(1-x)\left[(1-a-b)(c-1)+ab\right]+(c-1)(1+a+b-c)-ab}{x(1-x)(c-2)}y_n-\Frac{(1-c)}{x(1-x)} y_{n-1}\\
\\
y'_{n-1}= &\Frac{1-a-b+ab}{(1-x)(c-2)} y_{n-1}
-\Frac{x(a-c+1)(1-a)(c-b-1)(b-1)}{(1-x)(c-1)(c-2)^2}y_n\\
\end{array}
$$

Change of variable: $z(x)=-\displaystyle\sqrt{\Frac{(c-a-1)(1-a)(1+b-c)(b-1)}{(c-2)^2}}\ln (1-x)$.

\item{$DDE(1,0,1)$ :}
$$
\begin{array}{ll}
y'_n=& \Frac{1+xb-c}{x(1-x)}y_n-
\Frac{(1-c)}{x(1-x)}
  y_{n-1}\\
\\
y'_{n-1}= & -\Frac{(1-a) }{(1-x) } y_{n-1}
+\Frac{(1-a)(b+1-c)}{(1-x)(1-c)}y_n\\
\end{array}
$$

Change of variable: $z(x)=2\sqrt{(1-a)(b+1-c)}\tanh^{-1} \sqrt{x}$.

\item{$DDE(1,-1,0)$ :}
$$
\begin{array}{ll}
y'_n=&b\Frac{x(b-a+1)+a-c}{x(1-x)(b-a+1)}y_n+
\Frac{b(c-a)}{x(1-x)(b-a+1)}y_{n-1}\\
\\
y'_{n-1}= & (1-a)\Frac{(1-x)(b-a+1)+a-c}{x(1-x)(b-a+1)} y_{n-1}
-\Frac{(1-a)(1+b-c)}{x(1-x)(b-a+1)}y_n\\
\end{array}
$$

Change of variable: $z(x)=
\Frac{\sqrt{b(c-a)(1-a)(1+b-c)}}{b-a+1}\ln\left(\Frac{x}{1-x}\right)$.

\item{$DDE(0,0,-1)$ :}
$$
\begin{array}{ll}
y'_n=&\Frac{b+a-c}{1-x}y_n-\Frac{(b-c)(c-a)}{(1-x)c }y_{n-1}\\
\\
y'_{n-1}= & -\Frac{c}{x} y_{n-1}+\Frac{c}{x}y_n\\
\end{array}
$$

Change of variable: $z(x)=\sqrt{(b-c)(c-a)} \arcsin (2x-1)$.

\item{$DDE(1,1,1)$ :}
$$
\begin{array}{ll}
y'_n=&\Frac{x(a+b-1)-c+1}{x(1-x)}y_n-\Frac{1-c}{x(1-x)}y_{n-1}\\
\\
y'_{n-1}= & \Frac{(b-1)(1-a)}{(1-c)}y_n\\
\end{array}
$$

Change of variable: $z(x)=\sqrt{(b-1)(1-a)} \arcsin (2x-1)$.

\end{enumerate}

\subsubsection{Comparison of FPIs}

   The following table shows the $|D_{n}|=|-d_n e_n|$ coefficient for the FPIs
$(1,0,0)$, $(1,1,0)$, $(1,1,2)$, $(1,0,1)$, $(1,-1,0)$, $(0,0,-1)$ and $(1,1,1)$.
$$
\begin{array}{l|c}
\mbox{Iteration} & |D_n|  \\
\hline
 (1,0,0) &   \displaystyle\left|\Frac{(a-c)(a-1)}{x^2(1-x)}\right|\\
\hline
(1,1,0) &    \displaystyle\left|\Frac{(b-c)(a-c)(b-1)(a-1)}{x^2(a+b-c-1)^2}\right|        \\
\hline
(1,1,2) &    \displaystyle\left|\Frac{(c-a-1)(1-a)(1+b-c)(b-1)}{(1-x)^2(c-2)^2}\right|        \\
\hline
(1,0,1) &    \displaystyle\left|\Frac{(a-1)(1+b-c)}{x(1-x)^2}\right|     \\
\hline
(1,-1,0) &    \displaystyle\left|\Frac{b(c-a)(a-1)(c-b-1)}{x^2(1-x)^2(a-b-1)^2}\right|        \\
\hline
(0,0,-1) &     \displaystyle\left|\Frac{(b-c)(a-c)}{x(1-x)}\right|      \\
\hline
(1,1,1)  &      \displaystyle\left|\Frac{(b-1)(a-1)}{x(1-x)}\right|
\end{array}
$$

As can be inferred from the table, the most appropriate FPIs in
the interval $(0,1)$ are the $(0,0,-1)$ and $(1,1,1)$ iterations.
As commented, some hypergeometric functions in this interval with
particular values of their parameters are orthogonal polynomials
(Jacobi and derived polynomials).
In the case of orthogonal polynomials, the ``natural'' iteration to be considered is $(1,-1,0)$
(or $(-1,1,0)$ equivalently) which is not the optimal iteration. In order to illustrate this fact,
let us consider the evaluation of the zeros of the hypergeometric function
$_2\mbox{F}_1(-50,54;5/2;x)$ in the interval $(0,1)$. The zeros of this function
correspond to the zeros of the Jacobi polynomial $P^{(3/2,3/2)}_{50}(1-2x)$.
 In Figure 5, we show the ratio between the number of iterations
needed by FP$(1,-1,0)$ and FP$(1,1,1)$ as a function of the location of the zeros.

On the contrary, the most appropriate iteration in the $(1,
\infty)$ interval is the $(1,0,0)$ iteration. This could also be
understood taking into account that the best iterations in the
interval $(0,1)$ are FPI(0,0,-1) and FPI(1,1,1) and using Eq.\
(\ref{x>1}). In Figure 6  we show the ratio between the number of
iterations needed by FP$(1,1,1)$ and FP$(1,0,0)$.

The main conclusion for $_2\mbox{F}_1$ hypergeometric functions is
that the iteration $(1,1,1)$ is the preferred one and that
$(0,0,-1)$ can be considered as a replacement, with similar
performance. For the other two intervals, the relations
(\ref{x<0}) and (\ref{x>1}) indicate that the most appropriate
iterations will be $(1,0,1)$ for $(-\infty, 0)$ and $(1,0,0)$ for
$(1,+\infty)$; in any case, the solutions in these intervals can
be related to solutions in $(0,1)$.

\vspace*{0.2cm}

\centerline{\protect\hbox{\psfig{file=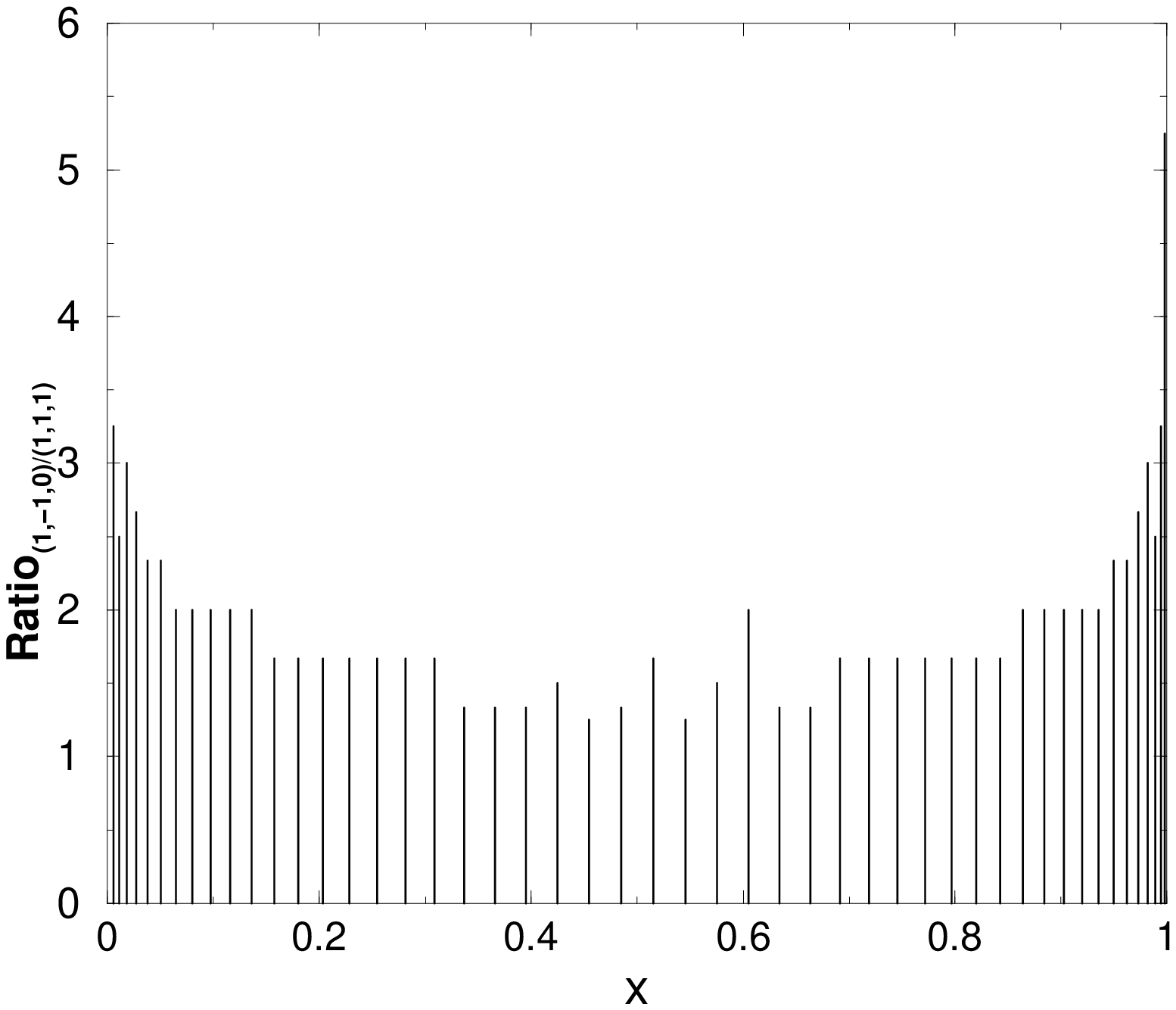,width=7.cm,height=6.cm}}}

{{\bf Figure 5}. Ratio between the number of iterations
needed by FP(1,-1,0) and FP(1,1,1) for the calculation of the zeros of the hypergeometric
function $_2\mbox{F}_1(-50,54;5/2;x)$, as a function of the location of the zeros.}

\centerline{\protect\hbox{\psfig{file=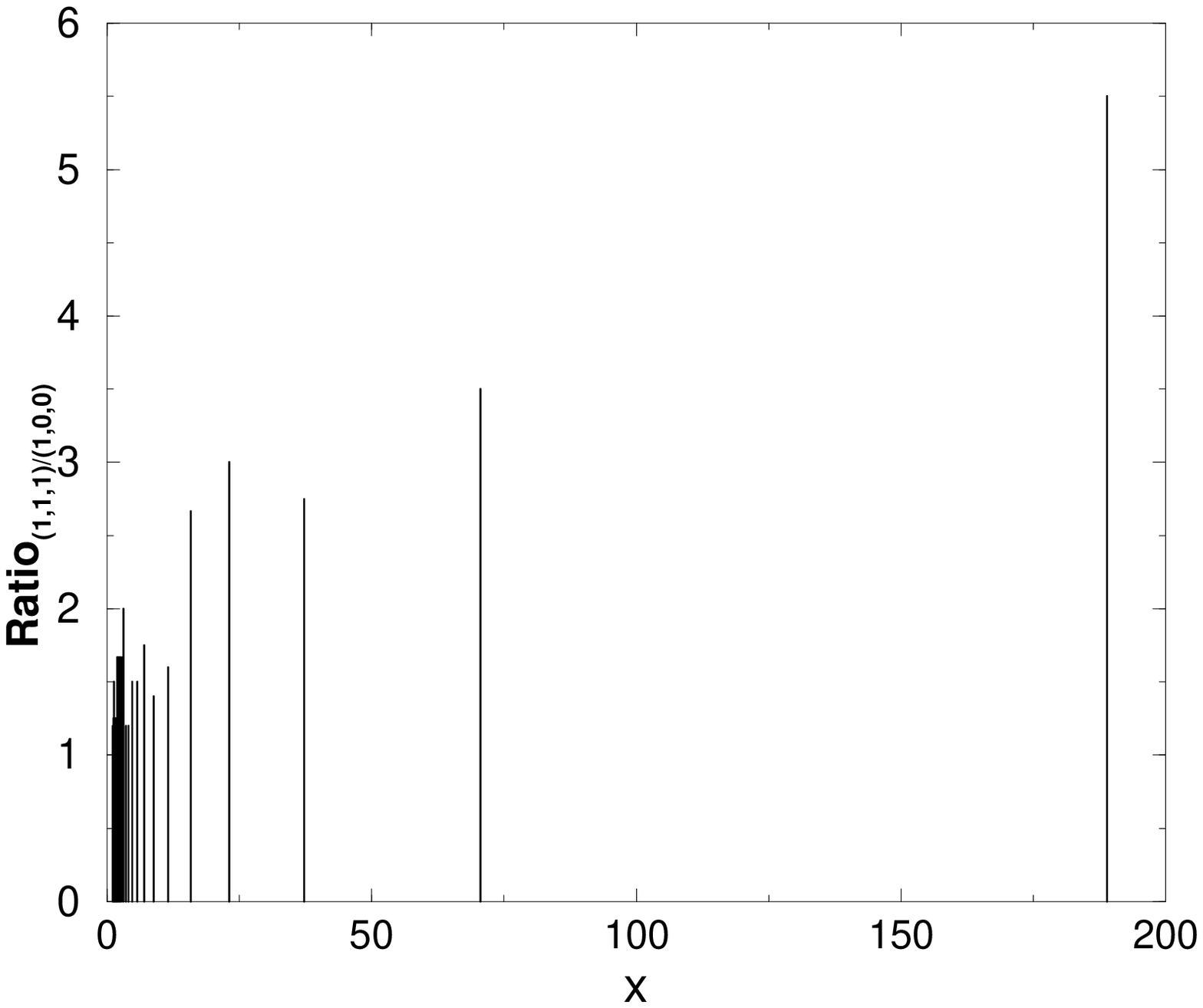,width=7.cm,height=6.cm}}}

{{\bf Figure 6}. Ratio between the number of iterations
needed by FP(1,1,1) and FP(1,0,0) for the calculation of the zeros of the hypergeometric
function $_2\mbox{F}_1(-30,-32;-70;x)$, as a function of the location of the zeros.}

\section{Conclusions}

We have developed a detailed study of the performance of the available fixed point methods for hypergeometric
 functions. This will allow the construction of efficient algorithms for the computation of the real
zeros of hypergeometric functions with a good asymptotic behavior.
The next table summarized the main results.

$$
\begin{array}{cccc}
\mbox{Function} & \mbox{Contrast function} & \mbox{Range of application} &
\mbox{Change of variables } z(x)\\
\hline
_0\mbox{F}_1 (;c;-x) & _0\mbox{F}_1 (;c_{-};-x) &
\left\{\begin{array}{l}
c<100\\
c>100\,,\,x>c^2/2
\end{array}
\right. & 2\sqrt{x}\\
`` & _0\mbox{F}_1 (;c_{-};-x) & c>100\,,\,x<c^2 /2 & x/(\nu-1)\\
\hline
_1\mbox{F}_1(a;c;x) & _1\mbox{F}_1 (a_{-};c_{-};x) & x<c-a & 2\sqrt{(1-a)}x\\
`` & _1\mbox{F}_1 (a_{-};c;x) & x>c-a & N_{ac} \ln x \\
\hline
_2\mbox{F}_1 (a,b;c;x) & _2\mbox{F}_1 (a_{-},b_{-};c_{-};x) & & N_{bc} \arcsin(2x-1)\\
\hline
\end{array}
$$
\noindent where $N_{ac}=\sqrt{(c-a)(1-a)}$ and $N_{bc}=\sqrt{(b-1)(1-a)}$ and
$a_{-}=a-1$, $b_{-}=b-1$, $c_{-}=c-1$.

In the table we restrict ourselves to $x>0$ in all cases, with the
additional restriction $x<1$ for the $_2\mbox{F}_1$ functions. For
the rest of the intervals, as commented before, relations are
available which map these other regions into the intervals
considered in the table.

It should be noted that, when a priori approximations to the roots
are available (for instance, asymptotic approximations like in
\cite{Tem}) the performance of the algorithms can be improved.
However, the methods presented here have the advantage of being
efficient methods that do not require specific approximations for
specific functions (which, on the other hand, are difficult to
obtain for three parameter functions like the $_2\mbox{F}_1$
hypergeometric functions). In addition, even in the simple cases
of one parameter functions, the methods are very efficient by
themselves.

To conclude, it is worth mentioning that one on the main reasons
for the good performance of the algorithms if that the analytical
transformations of the DDEs, and in particular, the associated
change of variable $z(x)=\int\sqrt{-d_n e_n}dx $, tend to
uniformize the distance between zeros. Generally speaking, the
most successful  methods are those which produce smaller
variations of the distances between zeros, because the first
guesses for the zeros become more accurate. In connection to this,
these changes of variable lead to interesting analytical
information about these zeros \cite{Dea}.

\begin{appendix}{Maple Code}
In this appendix we would like to  explain how we received the DDEs
automatically. Our Maple code uses a Zeilberger type approach
\cite{Zei} and is completely on the lines of \cite{Koe}. The
authors provide a Maple program {\it rules.mpl} which can be used
in combination with {\it hsum.mpl} \cite{Koe} to get equations
(15), (18), (22), (24), (26), (28), and the DDEs in \S~3.4.1.
These computations are collected in the Maple worksheet {\it
rules.mpl}. All these files can be obtained from the web site
\url{http://www.mathematik.uni-kassel.de/~koepf/Publikationen}.
\end{appendix}

\begin{acknowledgment}
  A. Gil acknowledges support from A. von Humboldt foundation. J. Segura acknowledges support from DAAD.
\end{acknowledgment}

\pagebreak

\end{article}
\end{document}